\documentclass[12pt]{article}
\usepackage{amssymb}
\usepackage{latexsym,bm}
\usepackage{graphicx}
\usepackage{amsmath}
\usepackage{mathrsfs}
\usepackage{epstopdf}

\date{}
\newcounter{mathitem}
  {\begin{list}{{$(\roman{mathitem})$}}{
   \setcounter{mathitem}{0}
   \usecounter{mathitem}
   \setlength{\topsep}{0pt plus 2pt minus 0pt}
   \setlength{\parskip}{0pt plus 2pt minus 0pt}
   \setlength{\partopsep}{0pt plus 2pt minus 0pt}
   \setlength{\parsep}{0pt plus 2pt minus 0pt}
   \setlength{\leftmargin}{35pt}
   \setlength{\itemsep}{0pt plus 2pt minus 0pt}}}
  {\end{list}}

\begin{document}
\title{The minimum size of a $k$-connected locally nonforesty graph
\footnote{E-mail addresses: {\tt lichengli0130@126.com} (C. Li), {\tt tyr2290@163.com} (Y. Tang), {\tt zhan@math.ecnu.edu.cn} (X. Zhan).}}
\author{\hskip -10mm Chengli Li, Yurui Tang and Xingzhi Zhan\thanks{Corresponding author}\\
{\hskip -10mm \small Department of Mathematics,  Key Laboratory of MEA (Ministry of Education) }\\
{\hskip -10mm \small \& Shanghai Key Laboratory of PMMP, East China Normal University, Shanghai 200241, China}}\maketitle
\begin{abstract}
A local subgraph of a graph is the subgraph induced by the neighborhood of a vertex. Thus a graph of order $n$ has $n$ local subgraphs. A graph $G$ is called locally nonforesty
 if every local subgraph of $G$ contains a cycle. Clearly, a graph is locally nonforesty if and only if every vertex of the graph is the hub of a wheel. We determine the minimum size of a $k$-connected locally nonforesty graph of order $n.$
\end{abstract}


\vskip 8mm

\section{Introduction}

We consider finite simple graphs and use standard terminology and notation from [2] and [9]. The {\it order} of a graph is its number of vertices, and the
{\it size} its number of edges.

{\bf Definition 1.} A {\it local subgraph} of a graph $G$ is the subgraph of $G$ induced by the neighborhood of a vertex.

Thus a graph of order $n$ has $n$ local subgraphs. It has been a traditional topic to deduce properties of a graph by its local subgraphs. A graph $G$ is called {\it locally connected} if all local subgraphs of $G$ are connected. {\it Locally hamiltonian} graphs are defined similarly [5]. A well-known theorem of Oberly and Sumner [7] states that every connected, locally connected claw-free graph is hamiltonian. Further connectivity conditions on such graphs force them to be hamiltonian-connected ([1] and [3]). An interesting unsolved conjecture of Ryj\'{a}\v{c}ek (see [8]) asserts that every locally connected graph is weakly pancyclic.

{\bf Definition 2.} A graph $G$ is called {\it locally foresty} if every local subgraph of $G$ is a forest. $G$ is called {\it locally nonforesty} if every local subgraph
of $G$ contains a cycle.

Clearly, a graph is locally nonforesty if and only if every vertex of the graph is the hub of a wheel. Recently, in studying forest cuts of a graph, Chernyshev, Rauch and Rautenbach [4, p.8] posed the following

{\bf Conjecture 1.} {\it If $n$ and $m$ are the order and size of a $3$-connected locally nonforesty graph respectively, then $m\ge 7(n-1)/3.$}

In [6] we have solved this problem by determining the minimum size of a $3$-conneted locally nonforesty graph of order $n.$ It turned out that Conjecture 1 does not hold.

In this paper we determine the minimum size $f(k,n)$ of a $k$-connected locally nonforesty graph of order $n$ for a general $k.$
If $k\ge 5,$ the problem is trivial: $f(k,n)=\lceil kn/2\rceil;$ i.e., the condition of being locally nonforesty has no effect on the minimum size. For $k=5$ we have constructed
$5$-connected locally nonforesty graphs of order $n$ and the size $\lceil kn/2\rceil,$ but it is tedious to describe those graphs and verify their properties. For $k\ge 6,$ Harary's graphs [9, pp.150-151] do the job. Hence we treat only the nontrivial cases $k\le 4$ of the problem.

We denote by $V(G)$ the vertex set of a graph $G,$ and denote by $e(G)$ the size of $G.$ The neighborhood of a vertex $x$ is denoted by $N(x)$ or $N_G(x),$ and the
 closed neighborhood of $x$ is $N[x]\triangleq N(x)\cup \{x\}.$ The degree of $x$ is denoted by ${\rm deg}(x).$ We denote by $\delta (G)$ and $\Delta(G)$ the minimum degree and maximum degree of $G,$ respectively.  For a vertex subset $S\subseteq V(G),$ we use $G[S]$ to denote the subgraph of $G$ induced by $S,$ and use $N(S)$ to denote the neighborhood
 of $S;$ i.e., $N(S)=\{y\in V(G)\setminus S \,|\, y\,\,{\rm has}\,\,{\rm a}\,\,{\rm neighbor}\,\,{\rm in}\,\,S\}.$ For $x\in V(G)$ and $S\subseteq V(G),$ $N_S(x)\triangleq N(x)\cap S$ and the degree of $x$ in $S$ is ${\rm deg}_S(x)\triangleq |N_S(x)|.$ Given two disjoint vertex subsets $S$ and $T$ of $G,$ we denote by $[S, T]$ the set of edges having one endpoint in $S$ and the other in $T.$  The degree of $S$ is ${\rm deg}(S)\triangleq |[S, \overline{S}]|,$ where $\overline{S}=V(G)\setminus S.$ We denote by $C_n,$ $P_n$ and $K_n$ the cycle of order $n,$ the path of order $n$ and the complete graph of order $n,$ respectively. $\overline{G}$ denotes the complement of a graph $G.$ For two graphs $G$ and $H,$ $G\vee H$ denotes the {\it join} of $G$ and $H,$ which is obtained from the disjoint union $G+H$ by adding edges joining every vertex of $G$ to every vertex of $H.$

For graphs we will use equality up to isomorphism, so $G=H$ means that $G$ and $H$ are isomorphic.

The paper is organized as follows. We treat $k$-connected graphs for $k=4,2,1$ in Sections $2, 3, 4,$ respectively. We use this order because our proof of the case $k=1$
relies on the result for $k=2.$

\section{The minimum size of a $4$-connected locally nonforesty graph}

{\bf Lemma 1.} {\it Let $G$ be a graph of order $n$ such that $G[N(v)]=C_3+K_1$ for every vertex $v\in V(G).$ Then $n\equiv 0~(\,{\rm mod}\,\, 4)$.}

{\bf Proof. } Clearly, $G$ is $4$-regular and each vertex in $G$ lies in a subgraph isomorphic to $K_4$. Let $A_1,\dots,A_t$ be all vertex sets in $G$ such that $G[A_i]=K_4$ for $i=1,2,\dots,t$. We claim that $A_j\cap A_k=\emptyset$ for any $j\neq k$. To the contrary, suppose $A_j\cap A_k\neq\emptyset$ for some $j\neq k$.. If $|A_j\cap A_k|=1,$ then there exists a vertex with degree at least $6$, a contradiction. If $|A_j\cap A_k|=2$, then there exists a vertex with degree at least $5$, a contradiction. If $|A_j\cap A_k|=3,$ then there exists a vertex $v\in A_j\cap A_k$ such that $G[N(v)]\neq C_3+ K_1$, a contradiction. Therefore, $n\equiv 0~(\,{\rm mod}\,\, 4)$.

{\bf Theorem 2.} {\it Given an integer $n\ge 8,$ define
	$$
	h(n)=\begin{cases}
		2n\quad\quad\quad\,{\rm if}\,\,\, n\equiv 0~(\,{\rm mod}\,\, 4),\\
		2n+1 \quad\,\,\,{\rm otherwise.}
	\end{cases}
	$$
	Then the minimum size of a $4$-connected locally nonforesty graph of order $n$ is $h(n).$
}

{\bf Proof. } Let $G$ be a $4$-connected locally nonforesty graph of order $n$. We first prove $e(G)\ge h(n).$ Since $\delta(G)\ge \kappa(G)\ge 4$, we have $e(G)\ge 2n.$
Hence to show $e(G)\ge h(n),$ it suffices to prove that $e(G)=2n$ occurs only when $n\equiv 0~(\,{\rm mod}\,\, 4)$.

Now suppose $e(G)=2n.$ Then $G$ is $4$-regular, since $\delta(G)\ge 4.$ For a vertex $v$ of $G,$ we denote by $L(v)$ the subgraph of $G$ induced by $N(v).$
We assert that $L(v)$ contains no $C_4$. To the contrary, suppose $L(v)$ contains a $C_4$. If $L(v)=K_4$, since $G$ is $4$-regular, then $G=K_5$, a contradiction. If $L(v)=K_4-xy$ where $xy$ is an edge, then $\{x, y\}$ is a vertex cut of $G$, contradicting the condition that $G$ is $4$-connected. Suppose $L(v)=C_4=v_1v_2v_3v_4v_1$. Let $N(v_1)=\{v, v_2, v_4,w\}.$ Note that $w$ is nonadjacent to $v$ since $v$ already has four neighbors and $v_2$ is nonadjacent to $v_4$ since $L(v)$ is the $4$-cycle $v_1v_2v_3v_4v_1.$ Thus, the condition that $L(v_1)$ contains a cycle implies that $w$ is adjacent to both $v_2$ and $v_4.$ But then $\{w, v_3\}$ is a vertex cut of $G,$ a contradiction. This shows that
$L(v)$ contains no $C_4$ and hence $L(v)$ contains $C_3$, say $u_1u_2u_3u_1$. Let $u_4$ be the fourth neighbor of $v$. If $u_4$ is adjacent to one of $u_1, u_2, u_3$, say $u_1,$
then $\{u_2, u_3, u_4\}$ is a vertex cut of $G,$ a contradiction. Therefore, $L(v)=C_3+K_1$. Since the vertex $v$ was chosen arbitrarily, by Lemma 1, $n\equiv 0~(\,{\rm mod}\,\, 4).$

Now we construct a $4$-connected locally nonforesty graph $G_n$ of order $n$ with size $h(n).$
Given an integer $k\ge 2,$ let $A_i$ be a graph isomorphic to $K_4$ with vertex set $V(A_i)=\{x_i, y_i, z_i, w_i\}$ for $i=1,\ldots, k.$ The four auxiliary graphs
$A_i,$ $B_1,$ $C_1$ and $D_1$ are depicted in Figure 1.
\begin{figure}[h]
\centering
\includegraphics[width=0.7\textwidth]{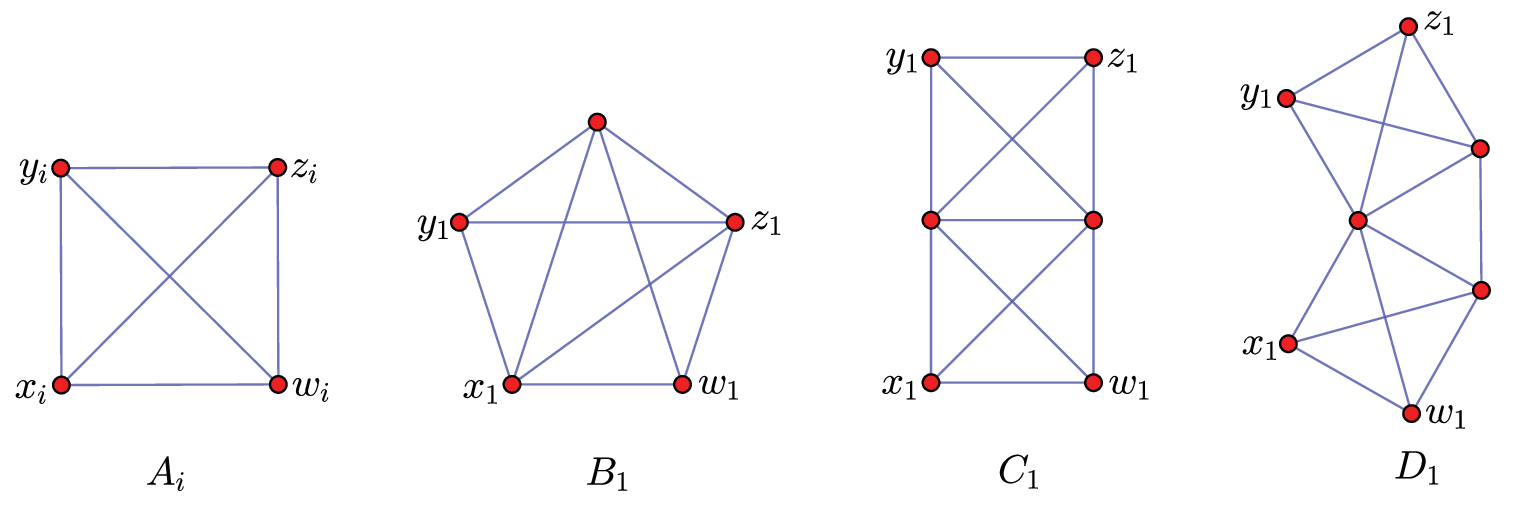}
\caption{$A_i,$ $B_1,$ $C_1$ and $D_1$}
\end{figure}
Suppose $n=4k.$ Let $G_n$ be the graph obtained from the disjoint union $A_1+A_2+\cdots +A_{k}$ by adding edges $z_iy_{i+1},w_ix_{i+1},$ $1\le i\le k$ where $y_{k+1}=y_1,x_{k+1}=x_1$.

Suppose $n=4k+1.$ Replace $A_1$ by $B_1$ and then construct $G_n$ as above.

Suppose $n=4k+2.$ Replace $A_1$ by $C_1$ and then construct $G_n$ as in the case $n=4k$.

Suppose $n=4k+3.$ Replace $A_1$ by $D_1$ and then construct $G_n$ as in the case $n=4k$.

\section{The minimum size of a $2$-connected locally nonforesty graph}

In this section we determine the minimum size of a $2$-connected locally nonforesty graph of order $n.$

{\bf Theorem 3.} {\it Given an integer $n\ge 8,$ define
$$
g(n)=\begin{cases}
2n-\lfloor n/4\rfloor\quad\quad\quad\,{\rm if}\,\,\, n\equiv 0,3~(\,{\rm mod}\,\, 4),\\
2n+1-\lfloor n/4\rfloor \quad\,\,\,{\rm otherwise.}
\end{cases}
$$
Then the minimum size of a $2$-connected locally nonforesty graph of order $n$ is $g(n).$
}

{\bf Proof.} Let $G$ be a $2$-connected locally nonforesty graph of order $n.$ We first prove that $e(G)\ge g(n).$
Let $n=4k+r$ with $0\le r\le 3.$ Then
$$
g(n)=\begin{cases}
7k+2r~\quad\quad\quad\quad\,{\rm if}\,\,\, r=0~{\rm or}~3\\
      7k+2r+1 ~\quad\quad\,\,\,{\rm if}\,\,\, r=1~{\rm or}~2.
\end{cases}
$$
For a vertex $v$ of $G,$ we denote by $L(v)$ the subgraph of $G$ induced by $N(v).$

Since $G$ is a locally nonforesty graph, $\delta(G)\ge 3.$ If $\delta(G)\ge 4,$ we have
$e(G)\ge 2n\ge g(n).$ Next we suppose $\delta(G)=3.$ Denote
$$
S=\{v\in V(G)|\,{\rm deg}(v)=3\},\quad T=N(S),\quad W=V(G)\setminus (S\cup T).
$$
We assert that $\Delta(G[S])\le 1$, i.e., $G[S]$ is a graph which is a union of some isolate vertices and a matching. To the contrary, suppose that $S$ contains a vertex $u$ with degree at least two in $G[S].$  Since $L(u)$ contains a cycle, $G[N[u]]=K_4$. Thus $u$ has exactly degree two in $G[S]$, since otherwise, $G=K_4$, a contradiction. Then the vertex in
 $N(u)\setminus N_S(u)$ would be a cut-vertex, contradicting our assumption that $G$ is $2$-connected.

Let $T_i=\{v\in T|\, {\rm deg}_S(v)=i\}$ where $1\le i\le s\triangleq |S|.$ For any vertex $u\in T_i,$ we have ${\rm deg}(u)\ge i+1,$ which can be deduced from the facts
that $G[S]$ is a forest and that $L(u)$ contains a cycle. We distinguish two cases.

{\bf Case 1.} There exists an $i$ such that $T_i$ contains a vertex $u$ with ${\rm deg}(u)=i+1.$

Let $A=N_S(u)$ where $|A|=i$ and $N(u)\setminus S=\{v\}.$ For each $p\in A,$ $G[N[p]]=K_4,$ implying that $N(p)\subseteq N[u]=A\cup\{u,v\}.$
Since $\Delta(G[S])\le 1$, we deduce that $|N(p)\cap A|=1$, $\{u,v\}\subseteq N(p)$ for any $p\in A$ and $G[A]=(|A|/2){K_2}.$ If there is a vertex not in $A\cup \{u,v\},$ then $v$ would be a cut vertex, contradicting the condition that $G$ is $2$-connected. Hence $V(G)=A\cup \{u,v\},$ implying that $G=K_2\vee (|A|/2){K_2}.$ It follows that
$e(G)=(5n/2)-4>2n+1-\lfloor n/4\rfloor,$ where we have used the assumption that $n\ge 8.$

{\bf Case 2.} For every $i$ and any vertex $u\in T_i,$ ${\rm deg}(u)\ge i+2.$

Recall that we have denoted $s=|S|.$ Note that $\sum_{i=1}^s |T_i|=|T|$ and
$$
\sum_{i=1}^s i\cdot |T_i|=\sum_{v\in T}{\rm deg}_S(v)=|[S,\, T]|=\sum_{u\in S}{\rm deg}_T(u)=\sum_{u\in S}({\rm deg}_G(u)-{\rm deg}_S(u))= 3s-2e(G[S])\ge 2s.
$$
Observe that every vertex in $T\cup W$ has degree at least $4.$ We have
$$
e(G)\ge (3s+4(n-s))/2=2n-s/2
$$
and
\begin{align*}
e(G)&\ge (3s+\sum_{i=1}^{s} (i+2)|T_i|+4|W|)/2\\
    &\ge (5s+2|T|+4|W|)/2\\
    &\ge n+3s/2.
\end{align*}
Denote $\varphi(s)={\rm max}\{\lceil 2n-s/2\rceil,\,\lceil n+3s/2\rceil \}.$ Then
$$
e(G)\ge \mathop{\rm min}\limits_{1\le s\le n}\varphi(s). \eqno (2)
$$

According to the value of the remainder $r$ in $n=4k+r$ with $0\le r\le 3,$ we distinguish four subcases.

{\bf Subcase 2.1.} $n=4k.$

By inequality (2), we have $e(G)\ge \mathop{\rm min}\limits_{1\le s\le n}\varphi(s)=\varphi(2k)=7k=g(n).$

{\bf Subcase 2.2.} $n=4k+1.$

In this case, $\mathop{\rm min}\limits_{1\le s\le n}\varphi(s)=7k+2$ and the minimum value $7k+2$ is attained uniquely at $s=2k.$ By (2) we have
$e(G)\ge 7k+2$ and equality holds only if $s=2k.$ Now we exclude the possibility that $e(G)=7k+2,$ and then it follows that $e(G)\ge 7k+3=g(n).$

To the contrary, suppose $e(G)=7k+2.$ Then $s=2k.$ It follows that the degree sequence of $G$ is
$(\operatorname*{\underbrace{3,\cdots, 3}}\limits_{2k},\operatorname*{\underbrace{4,\cdots, 4}}\limits_{2k+1}).$ Let $S=\{a_1,\ldots,a_{2k}\},$
the set of vertices of degree $3.$ Denote $A_i=N[a_i],$ the closed neighborhood of $a_i.$ Then $G[A_i]=K_4$ for each $i=1,\ldots,2k.$

{\bf Claim 1.} Either $A_i\cap A_j=\emptyset$ or $A_i=A_j$ for $i\neq j.$

To the contrary, suppose  $A_i\cap A_j\neq\emptyset$ and $A_i\neq A_j.$
If $|A_i\cap A_j|=1$ or $2,$ then  $A_i\cap A_j$ contains a vertex of degree at least $5,$ a contradiction. If $|A_i\cap A_j|=3$ then $n=5,$ contradicting our
assumption that $n\ge 8.$

Since $\Delta(G[S])\le 1$, each $A_i$ contains at most two vertices in $S$. Using Claim 1 and the condition that $n=4k+1,$ it is easy to show that among $A_1,\ldots, A_{2k}$ there are exactly $k$ sets that are pairwise distinct. Then each $A_i$ contains exactly two vertices in $S$, that is, for each $i$, there is exactly one integer $j\neq i$ such that $A_i=A_j.$
Let $V(G)\setminus (\bigcup_{i=1}^{2k}A_i)=\{x\}.$

{\bf Claim 2.} $|N(x)\cap A_i|\le 1$ for every $i$ with $1\le i\le 2k.$

Note that $|A_i|=4$ and $a_i$ and $x$ are nonadjacent for each $i.$ Since each $A_i$ contains exactly two vertices in $S$, we have $|N(x)\cap A_i|\le 2.$ Using the conditions $\Delta(G)=4$ and ${\rm deg}(a_i)=3,$ we deduce that
if $|N(x)\cap A_i|=2$ then $x$ is a cut-vertex of $G,$ contradicting
the condition that $G$ is $2$-connected. This shows that $|N(x)\cap A_i|\le 1.$

Since $\Delta(G)=4$, Claim 2 implies that $N(x)$ is an independent set and hence $L(x)$ does not contain a cycle, a contradiction. This shows that $e(G)\neq 7k+2.$

{\bf Subcase 2.3.} $n=4k+2.$

$g(n)=7k+5.$ Now $\mathop{\rm min}\limits_{1\le s\le n}\varphi(s)=7k+4$ and the minimum value $7k+4$ is attained at $s=2k$ or $s=2k+1.$ It suffices to
exclude the possibility that $e(G)=7k+4.$ To the contrary, suppose $e(G)=7k+4.$ Then $s=2k$ or $s=2k+1.$

(a) Assume $s=2k.$
Now the degree sequence of $G$ is $(\operatorname*{\underbrace{3,\cdots, 3}}\limits_{2k},\operatorname*{\underbrace{4,\cdots, 4}}\limits_{2k+2}).$ The proof in this case is similar to the one in Subcase~2.2, and we omit the details.

(b) Assume $s=2k+1.$
Now the degree sequence of $G$ is $(\operatorname*{\underbrace{3,\cdots, 3}}\limits_{2k+1},\operatorname*{\underbrace{4,\cdots, 4}}\limits_{2k},5).$
Denote $A_i=N[a_i],$ the closed neighborhood of $a_i.$ Then $G[A_i]=K_4$ for each $i=1,\ldots,2k+1.$

{\bf Claim 3.} Either $A_i\cap A_j=\emptyset$ or $A_i=A_j$ for $i\neq j.$

To the contrary, suppose  $A_i\cap A_j\neq\emptyset$ and $A_i\neq A_j.$
If $|A_i\cap A_j|=1$,  then  $A_i\cap A_j$ contains a vertex of degree at least $6,$ a contradiction. If $|A_i\cap A_j|=2$,  then  $A_i\cap A_j$ contains at least two vertices of degree at least $5,$ a contradiction. If $|A_i\cap A_j|=3$, since $n\ge 8,$ the vertex of degree five is a vertex-cut, a contradicting our assumption that $G$ is $2$-connected.

Note that each $A_i$ contains at most two vertices in $S$. Then there are at least $k+1$ disjoint subgraphs isomorphic to $K_4$, since $s=2k+1.$ Thus $n\ge 4k+4$, a contradiction. This shows that $e(G)\neq 7k+4.$

{\bf Subcase 2.4.} $n=4k+3.$

By inequality (2), we have $e(G)\ge \mathop{\rm min}\limits_{1\le s\le n}\varphi(s)=7k+6=g(n).$

Now for every integer $n\ge 8,$ we construct a $2$-connected locally nonforesty graph $G_n$ of order $n$ and size $g(n).$
Given an integer $k\ge 2,$ let $A_i$ be a graph isomorphic to $K_4$ with vertex set $V(A_i)=\{x_i, y_i, z_i, w_i\}$ for $i=1,\ldots, k.$ The graphs $B_1,$ $C_1$
and $D_2$ are depicted in Figure 2.
\begin{figure}[h]
\centering
\includegraphics[width=0.8\textwidth]{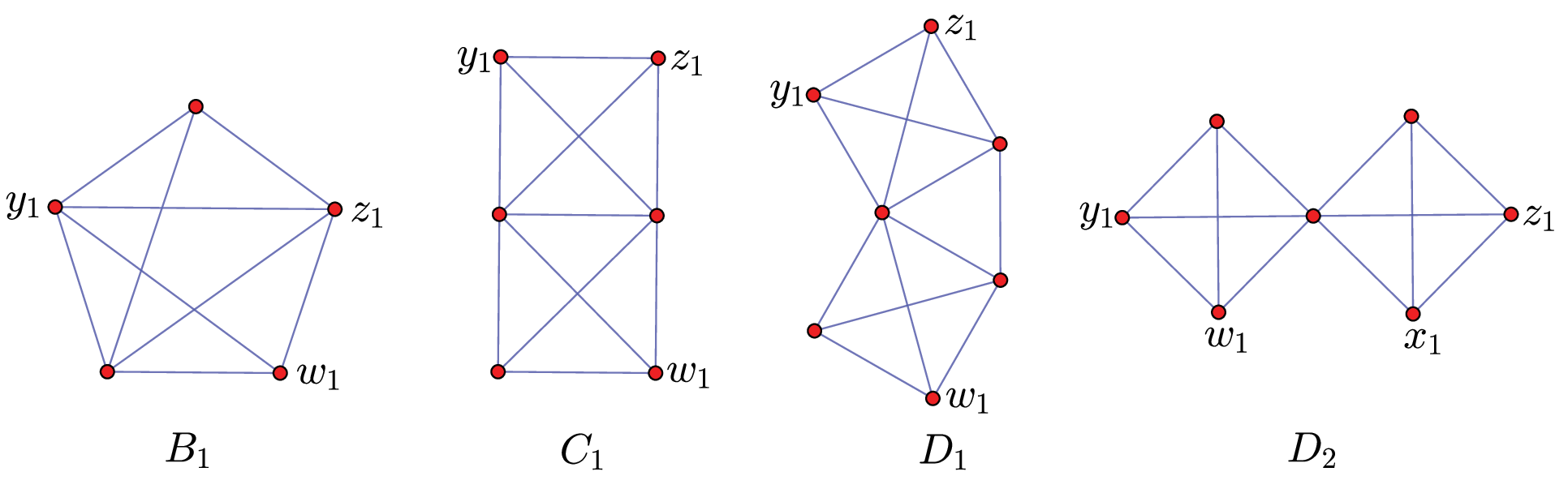}
\caption{$B_1,$ $C_1,$ $D_1$ and $D_2$}
\end{figure}

Suppose $n=4k.$ Let $G_n$ be the graph obtained from the disjoint union $A_1+A_2+\cdots +A_{k}$ by adding edges $z_iy_{i+1},$ $1\le i\le k$ where $y_{k+1}=y_1$.

Suppose $n=4k+1.$ Replace $A_1$ by $B_1$ and then construct $G_n$ as above.

Suppose $n=4k+2.$ Replace $A_1$ by $C_1$ and then construct $G_n$ as in the case $n=4k$.

Suppose $n=4k+3.$ Replace $A_1$ by $D_2$ and then construct $G_n$ as in the case $n=4k$.

\section{The minimum size of a connected locally nonforesty graph}

In this section we determine the minimum size of a connected locally nonforesty graph of order $n.$

{\bf Theorem 4.} {\it Given an integer $n\ge 8,$ define
$$
p(n)=\begin{cases}
2n-1-\lfloor n/4\rfloor ~~\quad\,{\rm if}\,\,\, n\equiv 0,3~(\,{\rm mod}\,\, 4),\\
2n-\lfloor n/4\rfloor \quad\quad\quad\,\,\,{\rm otherwise.}
\end{cases}
$$
Then the minimum size of a connected locally nonforesty graph of order $n$ is $p(n).$
}

{\bf Proof.}
Let $G$ be a connected locally nonforesty graph of order $n$ with minimum size. We first prove that $e(G)\ge p(n)$.

For a vertex $v$ of $G,$ we denote by $L(v)$ the subgraph of $G$ induced by $N(v).$ Let $H$ be a subgraph of $G$, we denote by $L_H(v)$ the subgraph of $H$ induced by $N_H(v).$

If $G$ is $2$-connected, then the result holds by Theorem 3.
Suppose that $G$ is not $2$-connected. Then $G$ has a block-cutpoint graph [9, p.156]. Let $M_1,M_2,\dots,M_t$ be all the blocks of $G$ of order at least $3.$ Let $m_i=|M_i|$ for $1\le i\le t$.

{\bf Claim 1.} $m_i\neq 3$ for each $i$.

To the contrary, suppose that there exists an integer $i$ such that $m_i=3$. Since $M_i$ is $2$-connected, $M_i=K_3$. For each vertex $u\in V(M_i)$, $L(u)$ contains a cycle that does not contain any edge of $M_i.$ Thus, by deleting an edge in $M_i$, we obtain a connected locally nonforesty graph with a smaller size, a contradiction.

{\bf Claim 2.} For each $i$, the graph $M_i$ is locally nonforesty.

To the contrary, suppose that there exists a vertex $v\in V(M_1)$ such that $L_{M_1}(v)$ is a forest.
Since $M_1$ is a block, we have $\deg_{M_1}(v)\ge 2$.
If $L_{M_1}(v)$ contains an isolated vertex, say $x$, then $N_G(x)\cap N_G(v)=\emptyset$ and hence $G-xv$ is a connected locally nonforesty graph with a smaller size, a contradiction.

Suppose that $L_{M_1}(v)$ contains an isolated edge, say $uu'$. Consider a new graph $G-uv$. The edge $uv$ does not lie in a cycle of $L(u')$, otherwise the degree of $u'$ in $L_{M_1}(v)$ is at least two, a contradiction. Note that $N_G(u)\cap N_G(v)=N_{M_1}(u)\cap N_{M_1}(v)=\{u'\}$. Then for each vertex $z\neq u'$, the edge $uv$ does not lie in a cycle of $L(z)$.
Therefore, $G-uv$ is a connected locally nonforesty graph with a smaller size, a contradiction.

Suppose that $L_{M_1}(v)$ does not contain an isolated vertex or an isolated edge.
Let $u$ be a leaf of $L_{M_1}(v)$ and let $u'$ be the unique neighbor of $u$ in $L_{M_1}(v)$. If $L(u')$ contains a cycle in $G-uv$, then $G-uv$ is a connected locally nonforesty graph with a smaller size, a contradiction. Suppose that $L(u')$ is a forst in $G-uv$. Let $C$ be a cycle of $L(u')$ in $G$ and let $u''$ be a neighbor of $v$ in $C$ that is distinct from $u.$ Since $L_{M_1}(v)$ is a forest, $u$ and $u''$ are nonadjacent. By deleting edge $uv$ and adding edge $uu''$, we obtain a new connected locally nonforesty graph $G'$ such that $e(G)=e(G')$ and $\deg_{G'}(v)=\deg_{G}(v)-1$. Clearly, the size of subgraph $L_{M_1}(v)$ is reduced by one in $G'$ compared to $G$.
Repeating the above operation, we obtain a connected locally nonforesty graph with the same size as $G$ and $L_{M_1}(v)$ contains an isolated edge. However, according to the above analysis, we also obtain a connected locally nonforesty graph with a smaller size, a contradiction.

By Claim 2, $M_i$ is a $2$-connected locally nonforesty graph. For $m_i\ge 8$, by Theorem 3, let $m_i=4k+r$ with $0\le r\le 3,$ we have
$$
e(M_i)\ge\begin{cases}
7k+2r~\quad\quad\quad\quad\,{\rm if}\,\,\, r=0~{\rm or}~3\\
      7k+2r+1 ~\quad\quad\,\,\,{\rm if}\,\,\, r=1~{\rm or}~2.
\end{cases}
$$
It is not difficult to check that $e(M_i)=6$ if $m_i=4$, $e(M_i)\ge 9$ if $m_i=5$, $e(M_i)\ge 11$ if $m_i=6$ and $e(M_i)\ge 13$ if $m_i=7$.

{\bf Claim 3.} For each $i$, we have $m_i\le 6$.

Ton the contrary, suppose there exists an integer $i$ such that $m_i\ge 7$. If $m_i=7$, then we can replace the block $M_i$ by $D_2$ which is depicted in Figure 2.
If $m_i\ge 8$, then we can replace $M_i$ by $G_n-z_1y_2$ where $G_n$ is defined in the proof of Theorem 3.
 Therefore, we obtain a connected locally nonforesty graph with a smaller size, a contradiction.

Let $t_i$ be the number of blocks whose order is $i$.

{\bf Claim 4.} $t_5+t_6\le 1$.

To the contrary, suppose that $t_5\ge 2.$ Other cases are similar. We can replace these two blocks of order 5 with two new blocks, one isomorphic to $K_4$ and the other isomorphic to $C_1$ which is depicted in Figure 4. Thus we obtain a connected locally nonforesty graph with a smaller size, a contradiction.

By Claims 1 and 3, we have that $t_3=0$ and $t_i=0$ for $i\ge 7$.
If $t_5=t_6=0$, then $e(G)=t_2+6t_4$ and $n=2t_2+4t_4-(t_2+t_4-1)=t_2+3t_4+1.$
Thus $e(G)=2n-2-t_2.$
Since $G$ is locally nonforesty, any two blocks of order 2 are nonadjacent and any block of order 2 is not an end block. Thus $t_2\leq t_4-1.$ Combine the fact that $n=t_2+3t_4+1$ and $t_2$ is an integer, we have $t_2\leq \lfloor(n-4)/4\rfloor$. Therefore,
$e(G)\geq 2n-1-\lfloor n/4\rfloor$.
Moreover, if $t_2\leq t_4-2$, then $t_2\leq \lfloor(n-7)/4\rfloor$ and hence $e(G)\geq 2n-1-\lfloor (n-3)/4\rfloor$.

{\bf Claim 5.} If $t_5=t_6=0$ and $t_2=t_4-1,$ we have $n\equiv 0~(\,{\rm mod}\,\, 4)$.

Note that $M_i=K_4$ for each $i=1,\dots,t.$ Then for every vertex $x$, $x$ lies in a subgraph isomorphic to $K_4$, since any two blocks of order 2 are nonadjacent
and any block of order 2 is not an end block. Since $t_2=t_4-1,$ we have $M_j\cap M_k=\emptyset$ for any $j\neq k$. Therefore, $n\equiv 0~(\,{\rm mod}\,\, 4)$.

If $t_5\neq 0$ or $t_6\neq 0$, by Claim 4, we have $t_5=1,t_6=0$ or $t_5=0,t_6=1.$
For the former case, we have $n=2t_2+4t_4+5-(t_2+t_4)=t_2+3t_4+5$
and $e(G)\ge t_2+6t_4+9=2n-1-t_2$.
Similarly as above, we have $t_2\leq t_4$ and hence
 $t_2\leq \lfloor(n-5)/4\rfloor$. Therefore,
$e(G)\geq 2n-\lfloor (n-1)/4\rfloor.$
For the latter case, we have $n=2t_2+4t_4+6-(t_2+t_4)=t_2+3t_4+6$
and $e(G)\ge t_2+6t_4+11=2n-1-t_2$.
Similarly as above, we have $t_2\leq t_4$ and hence
 $t_2\leq \lfloor(n-6)/4\rfloor$. Therefore,
$e(G)\geq 2n-\lfloor (n-2)/4\rfloor.$

Let $n=4k+r$ with $0\le r\le 3.$
According to the above, if $n=4k$ and $n=4k+3,$ then we have $e(G)\geq 2n-1-\lfloor n/4\rfloor$, as desired.

Suppose that $n=4k+1$. If $t_5=t_6=0$, by Claim~5, we have $t_2\leq t_4-2$ and hence $e(G)\geq 2n-1-\lfloor (n-3)/4\rfloor=2n-\lfloor n/4\rfloor.$ If $t_5=1,t_6=0$, then $e(G)\geq 2n-\lfloor (n-1)/4\rfloor=2n-\lfloor n/4\rfloor.$ If $t_5=0,t_6=1$, then $e(G)\geq 2n-\lfloor (n-2)/4\rfloor\ge 2n-\lfloor n/4\rfloor.$

Suppose that $n=4k+2$. If $t_5=t_6=0$, by Claim~5, we have $t_2\leq t_4-2$ and hence $e(G)\geq 2n-1-\lfloor (n-3)/4\rfloor=2n-\lfloor n/4\rfloor.$ If $t_5=1,t_6=0$, then $e(G)\geq 2n-\lfloor (n-1)/4\rfloor=2n-\lfloor n/4\rfloor.$ If $t_5=0,t_6=1$, then $e(G)\geq 2n-\lfloor (n-2)/4\rfloor=2n-\lfloor n/4\rfloor.$

Next for every integer $n\ge 8$, we construct a connected locally nonforesty graph $G_n$ of order $n$ and size $p(n).$ Thus $p(n)$ is indeed the minimum size.
Given an integer $k\ge 2,$ let $A_i$ be a graph isomorphic to $K_4$ with vertex set $V(A_i)=\{x_i, y_i, z_i, w_i\}$ for $i=1,\ldots, k.$ The graphs
$B_1,$ $C_1$ and $D_2$ are depicted in Figure 2.

Suppose $n=4k.$ Let $G_n$ be the graph obtained from the disjoint union $A_1+A_2+\cdots +A_{k}$ by adding edges $z_iy_{i+1},$ $1\le i\le k-1$.

Suppose $n=4k+1.$ Replace $A_1$ by $B_1$ and then construct $G_n$ as above.

Suppose $n=4k+2.$ Replace $A_1$ by $C_1$ and then construct $G_n$ as in the case $n=4k$.

Suppose $n=4k+3.$ Replace $A_1$ by $D_2$ and then construct $G_n$ as in the case $n=4k$.

{\bf Data availability}

No data was used for the research described in the article.

\vskip 5mm
{\bf Acknowledgement.} This research  was supported by the NSFC grant 12271170 and Science and Technology Commission of Shanghai Municipality
 grant 22DZ2229014.


\begin{thebibliography}{99}
\bibitem{1} A.S. Asratian, Every 3-connected, locally connected, claw-free graph is Hamilton-connected, J. Graph Theory, 23(1996), no.2, 191-201.
\bibitem{2} J.A. Bondy and U.S.R. Murty, Graph Theory, GTM 244, Springer, 2008.
\bibitem{3} G. Chartrand, R.J. Gould and A.D. Polimeni, A note on locally connected and Hamiltonian-connected graphs, Israel J. Math., 33(1979), no.1, 5-8.
\bibitem{4} V. Chernyshev, J. Rauch and D. Rautenbach, Forest cuts in sparse graphs, arXiv: 2409.17724, 26 September 2024.
\bibitem{5} J. Davies and C. Thomassen,  Locally Hamiltonian graphs and minimal size of maximal graphs on a surface, Electron. J. Combin., 27(2020), no.2, Paper No. 2.25.
\bibitem{6} C. Li, Y. Tang and X. Zhan, The minimum size of a $3$-connected locally nonforesty graph, arXiv:2410.23702, 31 October 2024.
\bibitem{7} D.J. Oberly and D.P. Sumner, Every connected, locally connected nontrivial graph with no induced claw is Hamiltonian, J. Graph Theory, 3(1979), no.4, 351-356.
\bibitem{8} D.B. West, Research problems, Discrete Math., 272(2003), 301-306.
\bibitem{9} D.B. West, Introduction to Graph Theory, Prentice Hall, Inc., 1996.
\end{thebibliography}
\end{document}